\documentstyle[12pt]{article}
\begin{document}

  {\large \bf Heegaard Diagrams of $S^3$ and the Andrews-Curtis Conjecture}

\vspace{0.6cm}

\hspace{2in}  {\large   Guangyuan Guo   }

\vspace{0.5cm}

 {\bf Abstract.}  {\em  We show that the Andrews-Curtis conjecture holds for
all balanced presentations of the trivial group corresponding to
Heegaard diagrams of $S^3$}

\vspace{0.5cm}


\vspace{0.5cm}


\vspace{1cm}

\hspace*{1.8in} {\large 1. \hspace{0.1cm} Introduction}

In this note, all manifolds are assumed to be orientable p.l.
manifolds and maps piecewise linear. If $G$ is a group and $g_{1},
\cdots, g_{k}$ are elements of $G$, then $(g_{1}, \cdots,
g_{k})$ and $\langle g_{1}, \cdots, g_{k} \rangle $ will denote
respectively the subgroup of $G$ generated by $g_{1},\cdots,
g_{k}$ and the smallest normal subgroup of $G$ containing $g_{1},
\cdots, g_{k}$. A set $R=\{r_{1}, \cdots, r_{n}\}$ of $n$ mutually
disjoint simple closed curves on a closed surface $S$ of
genus $n$ such that $S-\cup_{i} r_{i}$ is connected will be called
a {\bf complete system} on $S$ \cite{Ro}. A {\bf Heegaard diagram} is a
3-dimensional handlebody together with a complete system on its boundary
(our definition of a Heegaard diagram here is equivalent but slightly
different from that in some other papers, e.g. \cite{Ro}).

Let $n\geq 0$ be an integer and $V$ be a 3-dimensional
handlebody of genus $n$. {\bf A section of $V$} (also called a
complete system of meridian discs \cite{Pa} or a cut \cite{Zi}) is
a set $D=\{d_{1}, \cdots, d_{n}\}$ of $n$ mutually disjoint discs
in $V$ such that $\partial d_{i}=d_{i}\cap
\partial V$, $i=1, \cdots, n$, and that $V-\cup_{i}d_{i}$ is connected.
It is known that any section of $V$ can be obtained from any other
one by a finite number of operations which we call {\bf
replacements} (definition in section 2) in this note. This fact
gives a coarse relation between different Heegaard diagrams
associated with a Heegaard splitting, since each such diagram
corresponds to a section of a handlebody of the splitting.  We
observe that the operations above can be reduced to more
elementary ones called {\bf elementary replacements}. This gives a
more precise and tractable way to relate different Heegaard
diagrams associated with any given Heegaard splitting of a closed
orientable 3-manifold. It also enables us to define some equivalence
relation in the set of all complete systems on the boundary of a handlebody.
The set of equivalence classes of Heegaard splittings of genus $g$ are then
in one-to-one correspondence with the set of equivalence classes
of complete systems on the boundary of a handlebody of genus $g$.
And in particular, it provides a way to study all Heegaard
diagrams of $S^3$, since there is only one Heegaard splitting of
$S^3$ of genus $n$ up to equivalence, by a theorem of Waldhausen
\cite{Wa}. By examining the effect of each elementary replacement
on the elements of $\pi_{1}(V)$ represented by the boundary curves
of the discs in the sections of the complimentary handlebody in
$S^3$, we show that the Andrews-Curtis conjecture holds for all
balanced presentations of the trivial group corresponding to Heegaard
diagrams of $S^3$.

{\bf Acknowledgement} \hspace{2mm} The work in this note
was done while the author was holding positions at Nanchang University 
and the Central University of Finance and Economics. I thank both Universities for their
hospitality and generous financial support.

\vspace{1cm}

\hspace{4cm} {\large  2. \hspace{0.1cm} Sections of a Handlebody}

Let $V$ be a 3-dimensional handlebody of genus $n$ and $D=\{d_{1},
\cdots, d_{k}\}$ be a set of $k$ mutually disjoint discs in $V$
such that $\partial d_{i}=d_{i}\cap \partial V$, $i=1, \cdots, k$.
By {\bf cutting $V$ along $D$} (see \cite{He}) we mean removing
the interior of a regular neighborhood $N(d_{i})$ of each $d_{i}$,
$i=1, \cdots, k$, from $S$, where the $N(d_{i})$'s are always
chosen to be mutually disjoint. We will denote the result of
cutting $V$ along $D$ by $V_{D}$. $V_{\{d_{1}\}}$ will simply be
denoted by $V_{d_{1}}$. $V_{D}$ is a subset of $V$ and closed
subsets of $V$ that do not intersect the $d_{i}$'s can also be
considered as subsets of $V_{D}$. If $k\leq n$ and $V-
\cup_{i}d_{i}$ is connected, then $V_{D}$ is a handlebody of genus $n-k$. For each $i=1,
\cdots, k$, there are two obvious copies of $d_{i}$ on the
boundary of $ V_{D}$. Denote one copy by $d^{+}$, and the other by
$d^{-}$.

Given a section $D=\{d_{1}, \cdots, d_{n}\}$ of $V$, the set of
simple closed curves $\{\partial d_{1}, \cdots, \partial d_{n}\}$
on $\partial V$ is called the {\bf trace} of the section $D$
\cite{Zi}. Clearly, the trace of a section of $V$ is a complete
system on $\partial V$.

Now let $D=\{d_{1}, \cdots, d_{n}\}$ be a section of $V$ and let
$d_{1}'$ be a disc in $V$ such that $\partial d_{1}'=\partial V
\cap d_{1}'$, $d_{1}'\cap d_{i}=\emptyset$, $i=1, \cdots, n$, and
$V-d_{1}'$ is connected. On the surface of $V_{D}$, there are n
pairs $d_{i}^{\pm}$, $i=1, \cdots, n$, of discs corresponding to
the discs in $D$. Since $V_{D}$ is a 3-cell, it is decomposed into
two parts, each of which is a 3-cell, by $d_{1}'$. So
$(V_{D})_{d_{1}'}$ is the union of two disjoint 3-cells. If for
some $i$, $d_{i}^{+}$ and $d_{i}^{-1}$ are in separate parts, then
we say that $d_{1}'$ separates the pair $d^{\pm}_{i}$. Since
$V-d_{1}'$ is connected, $d_{1}'$ must separate some pair
$d_{i}^{\pm}$, $d_{1}^{\pm}$ say, and we can replace $d_{1}$ in
$D$ by $d_{1}'$ to get a new section $D'=\{d_{1}', d_{2}, \cdots,
d_{n}\}$ of $V$. We call such a process of getting a section $D'$
of $V$ from a given one $D$ a {\bf replacement}. If moreover,
$d_{1}'$ separates exactly two pairs and there are only two of the
$d_{i}^{\pm}$'s, $i=1, \cdots, n$, in one part of
$(V_{D})_{d_{1}'}$, then the replacement is called an {\bf
elementary replacement}.

{\bf Definition 2.1.} {\em Two sections $D_{1}$ and $D_{2}$ of a
handlebody $V$ are said to be equivalent, and written $D_{1}\sim
D_{2}$, if $D_{2}$ can be obtained from $D_{1}$ by a finite number
of elementary replacements.}

$``\sim"$ is obviously an equivalence relation in the set of all
sections of $V$. So the set of all sections of $V$ are divided
into equivalence classes. The following result says that there is
only one equivalence class of sections for any given handlebody.

{\bf Theorem 2.2.} {\em Let $V$ be a 3-dimensional handlebody and
$D=\{d_{1},\cdots, d_{n}\} $ and $D'=\{d_{1}', \cdots, d_{n}'\}$
be two sections of $V$. Then $D'$ can be obtained from $D$ by a
finite number of elementary replacements. }

{\bf Proof.} Theorem 2.2 is a refinement of Theorem 2.3 next,
using Lemma 2.4 that follows.

{\bf Theorem 2.3.} {\em Let $V$, $D$ and $D'$ be as in Theorem
3.2. Then $D'$ can be obtained from $D$ by a finite number of
replacements.}

{\bf Proof. } This is a known result. See, for example, Theorem 1
in \cite{Wj}.

{\bf Lemma 2.4.} {\em Let $V$, $D=\{d_{1},\cdots, d_{n}\}$ be as
above. Let $D'=\{d_{1}', d_{2}, \cdots, d_{n}\}$ be a section of
$V$ obtained from $D$ by one replacement (replacing $d_{1}$ by
$d_{1}'$). Then $D'$ can be obtained from $D$ by at most $n-1$
elementary replacements. }

{\bf Proof.} $d_{1}'$ separates $V_{D}$ into two parts, and by the
assumption, $d^{+}_{1}$ and $ d^{-}_{1}$ lie in separate parts.
Take a part that contains no more $d_{i}^{\pm}$'s, $i=1, \cdots,
n$, than the other part does and denote its closure by $P_{1}$.
$P_{1}$ contains exactly one of $d_{1}^{\pm}$, $d_{1}^{+}$ say,
and a copy of $d_{1}'$ which we still denote by $d_{1}'$. If
$P_{1}$ contains no other $d_{i}^{\pm}$, $2\leq i\leq n$, then
$d_{1}'$ is isotopic to $d_{1}$ in $V$. In this case $D$ and $D'$
are the same, up to isotopy, and no elementary replacement is
needed. If $P_{1}$ contains exactly one other $d_{i}^{\pm}, 2\leq
i\leq n$, then the replacement of $d_{1}$ by $d_{1}'$ is
elementary.

In general, suppose that $P_{1}$ contains a total of $m$
$d_{i}^{\pm}$'s, $i=1, \cdots, n$. By our choice of $P_{1}$,
$1\leq m\leq n$. Choose $m$ mutually disjoint discs
$\bar{d}_{1}=d_{1}^{+}, \bar{d}_{2}, \cdots, \bar{d}_{m-1},
\bar{d}_{m}=d_{1}'$ in $P_{1}$ such that $\partial\bar{d}_{j}=
\bar{d}_{j}\cap \partial P_{1}, \bar{d}_{j}\cap
d_{i}^{\pm}=\emptyset, 1\leq i\leq n, 2\leq j\leq m-1 $, and that
for each $j$, $2\leq j\leq m-1$, $\bar{d}_{1}, \cdots,
\bar{d}_{j-1}$ lie in the same side of $\bar{d}_{j}$ in $P_{1}$
($\bar{d}_{j}$ decomposes $P_{1}$ into two parts) and the side
contains exactly $j$ $d_{i}^{\pm}$'s. $\bar{d}_{2}, \cdots,
\bar{d}_{m-1}$ are also discs in $V$. Identify $\bar{d}_{1}$ with
$d_{1}$ and $\bar{d}_{m}$ with $d_{1}'$ in $V$. Then all
$\bar{d}_{1}, \cdots, \bar{d}_{m}$ are discs in $V$. Let $\bar
{D}_{j}=\{ \bar{d}_{j}, d_{2}, \cdots, d_{n}\}$, $ j=1, \cdots,
m$. Each $\bar{D}_{j}$ is a section of $V$, $\bar{D}_{1}=D,
\bar{D}_{m}=D'$ and $\bar{D}_{j}$, $j=2, \cdots, m$, is obtained
from $\bar{D}_{j-1}$ by replacing $\bar{d}_{j-1} $ by
$\bar{d}_{j}$ which is an elementary replacement. Thus we obtain
$D'=\bar{D}_{m}$ from $D=\bar{D}_{1}$ by $m-1\leq n-1$ elementary
replacements. This completes the proof of Lemma 2.4 and also the
proof of Theorem 2.2.

If $(U, V)$ is a Heegaard splitting of a closed 3-manifold $M$ and
$D=\{d_{1}, \cdots, d_{n}\}$ is a section of $U$, then $(R, V)$,
where $R=\{\partial d_{1}, \cdots,\partial d_{n}\}$, is a Heegaard
diagram for $M$. $(R, V)$ is said to be a Heegaard diagram
associated with the Heegaard splitting $(U, V)$. The set of all
Heegaard diagrams associated with the splitting $(U, V)$
corresponds to the set of all sections of $U$. By Theorem 2.2, any
section of $U$ can be obtained from any other section of $U$ by a
finite number of elementary replacements. Therefore any two
Heegaard diagrams associated with the splitting $(U, V)$ are
related through their corresponding sections. This gives a way to
study the relation between different Heegaard diagrams of a
Heegaard splitting.

Consider the set ${\cal C}$ of all complete systems on the
boundary $\partial V$ of the handlebody $V$. Define two types of
operations on elements of ${\cal C}$:

(1) $R=\{r_{1}, \cdots,r_{i-1}, r_{i}, r_{i+1}, \cdots, r_{n}\}
\rightarrow \tilde{R}=\{r_{1}, \cdots, r_{i-1}, \tilde{r}_{i},
r_{i+1}, \cdots, r_{n}\}$, where $\tilde{r}_{i}$ is a simple
closed curve obtained by sliding an $r_{j}$, $i\neq j$, over
$r_{i}$ along a simple curve $c$ joining $r_{i}$ and $r_{j}$ such
that the interior of $c$ does not meet any of the $r_{i}$'s.

(2) $R=\{r_{1}, \cdots, r_{n}\} \rightarrow h(R)=\{h(r_{1}),
\cdots, h(r_{n})\}$, where $h$ is a homeomorphism of $V$ onto
itself;

Now define an equivalence relation in ${\cal C}$ as follows. Two
elements $R$ and $\tilde{R}$ of ${\cal C}$ are said to be
equivalent if $R$ can be obtained from $R'$ by a finite number of
operations of types (1) and (2) above. One can easily check that
this is indeed an equivalence relation in ${\cal C}$. Since an
elementary replacement on the sections of a handlebody corresponds
to an operation of type (1) on the traces of the sections, it is
easy to see that Theorem 2.2 implies

{\bf Theorem 2.5. } {\em Two Heegaard diagrams $(R, V)$ and $(R',
V)$ are associated with equivalent Heegaard splittings if and only
if $R$ and $R'$ are equivalent in ${\cal C}$. Consequently, the
equivalence classes of Heegaard splittings of genus $n$ of closed
3-manifolds are in one-to-one correspondence with the equivalence
classes of complete systems on the boundary of a handlebody of
genus $n$.}

\vspace{1cm}

{\large 3. \hspace{0.1cm} Heegaard Diagrams of $S^3$ and the
Andrews-Curtis Conjecture}

By a theorem of Waldhausen \cite{Wa}, for each integer $n\geq 0$
there is only one Heegaard splitting for $S^3$ of genus $n$ up to
equivalence, therefore the result in the last section gives in
particular an effective way to investigate all Heegaard diagrams
of $S^3$. We look at this more closely.

For each $n\geq 0$, there is a canonical Heegaard diagram $(C,
W)$, $C=\{c_{1}, \cdots, c_{n}\}$, of genus $n$ for the 3-sphere
$S^{3}$. The diagram is characterized by the property that there
is a section $E=\{ e_{1}, \cdots, e_{n}\}$ of $W$ such that
$c_{i}\cap
\partial e_{j}$, $1\leq i, j\leq n$, is exactly one point if $i=j$
and empty if $i\neq j$. Let $(R, V)$ be any Heegaard diagram of
$S^3$. Consider the Heegaard splittings $(cl(S^{3}-V), V)$ and
$(cl(S^{3}-W), W)$ of $S^{3}$ associated with the diagrams $(R,
V)$ and $(C, W)$. By Waldhausen's theorem \cite{Wa}, there is a
homeomorphism $h:S^{3}\rightarrow  S^{3}$ such that $h(V)=W$ and
$h( cl(S^{3}-V))=cl(S^{3}-W)$. Using $h$ we identify $V$ with $W$,
which we denote by $V$, and $cl(S^{3}-V)$ with $cl(S^{3}-W)$,
which we denote by $U$. Then there are two sections $D_{R}$ and
$D_{C}$ of $U$, corresponding to $R$ and $C$ on $U\cap V$
respectively. By Theorem 2.2, $D_{R}$ can be obtained from $D_{C}$
by a finite number of elementary replacements. Therefore by
starting from the canonical diagram $(C, V)$ and studying the
effect of each elementary replacement on the traces of the
sections, one can hope to understand all Heegaard diagrams of
$S^3$.

{\bf Lemma 3.1.} {\em  Let $(U, V)$ be a Heegaard splitting of a
closed 3-manifold. Let $D'$ be a section of $U$ obtained from
another section $D$ of $U$ by an elementary replacement, then the
set of elements of $\pi_{1}(V, p)$ represented by the trace of
$D'$ may be obtained from the set of elements of $\pi_{1}(V, p)$
represented by the trace of $D$ by a finite number of
Andrews-Curtis transformations (Nielsen transformations and
conjugations).}

{\bf Proof.} Let $A=\{a_{1}, \cdots, a_{n}\}$ be the set of
elements of $\pi_{1}(V, p)$ represented by the trace of $D$ and
$B=\{b_{1}, \cdots, b_{n}\}$ the trace of $D'$. Suppose that
$D=\{d_{1},\cdots, d_{n}\}$ and $D'$ is obtained from $D$ by
replacing $d_{1}$ in $D$ by $d_{1}'$. So $D'=\{d_{1}', d_{2},
\cdots, d_{n}\}$ and $d_{1}'$ separates exactly two of the pairs
$d_{i}^{\pm}$, $i=1, \cdots, n$, in $U_{D}$ and moreover there are
only two of the $d_{i}^{\pm}$'s, one from the pair $d_{1}^{\pm}$
and the other from other pairs, in one of the two parts of
$(U_{D})_{d_{1}'}$. Without loss of generality, assume that the
other one of the $d_{i}^{\pm}$'s is from the pair $d_{2}^{\pm}$.
So $\partial d'$ is a simple closed curve obtained by sliding
$\partial d_{2}$ over $\partial d_{1}$ along a suitable simple
curve on $\partial U=\partial V$ connecting $\partial d_{1}$ and
$\partial d_{2}$. Then any element of $\pi_{1}(V, p)$ represented
by $\partial d'$ is a conjugate of
$a_{1}^{\epsilon_{1}}ga_{2}^{\epsilon_{2}}g^{-1}$, where $g$ is an
element of $\pi_{1}(V, p)$ depending on the curve used in sliding
and $\epsilon_{1}, \epsilon_{2}=1$ or $-1$ depending on the
orientation of $\partial d_{1}$, $\partial d_{2}$ and $\partial
d_{1}'$. In particular $b_{1}$ is a conjugate of
$a_{1}^{\epsilon_{1}}ga_{2}^{\epsilon_{2}}g^{-1}$. Consequently
$\{b_{1}, \cdots, b_{n}\}$ is conjugate to
$\{a_{1}^{\epsilon_{1}}ga_{2}^{\epsilon_{2}}g^{-1}, a_{2}, \cdots,
a_{n}\}$. Lemma 3.1 then follows.

{\bf Theorem 3.2.} {\em  If $(R, V)$ is a Heegaard diagram for
$S^{3}$ and $\{a_{1}, \cdots, a_{n}\}$ is a set of elements of
$\pi_{1}(V, p)$ represented by $R=\{r_{1}, \cdots, r_{n}\}$. Then
$\{a_{1}, \cdots, a_{n}\}$ can be transformed into a set of free
generators of the free group $\pi_{1}(V, p)$ by a finite number of
Andrews-Curtis transformations .}

{\bf Proof.} By the discussion preceding Lemma 3.1, $R$ is the
trace of a section $D_{R}$ of $U=cl(S^3-V)$. By Theorem 2.2,
$D_{R}$ can be transformed into the section $D_{C}$ by a finite
number of elementary replacements. It then follows from Lemma 4.1
that the set of elements $\{a_{1}, \cdots, a_{n}\}$ represented by
$R$ can be transformed into the set of elements of $F$ represented
by $C$ by a finite number of Andrews-Curtis transformations. The
set of elements of $\pi_{1}(V, p)$ represented by $C$ is, up to
conjugation, a set of free generators of $\pi_{1}(V, p)$. Thus
$\{a_{1}, \cdots, a_{n}\}$ can be transformed into a set of free
generators of $\pi_{1}(V, p)$ by a finite number of Andrews-Curtis
transformations. This completes the proof of Theorem 3.2.

Note that Theorem 3.2 does include some non-trivial case. By
trivial case we mean the case in which the set of elements of
$\pi_{1}(V, p)$ represented by $R$ is conjugate to a free basis.
The trivial case is precisely the case when the diagram $(R, V)$
is equivalent to the canonical diagram $(C, V)$. There are
diagrams for $S^3$ that are not equivalent to $(C, V)$.

Also note that in \cite{DR}, Rolfsen states that

{\bf Theorem:} {\em The AC conjecture is true for spines.}

What Rolfsen calls ``the AC conjecture" here is a geometric
version of the Andrews-Curtis conjecture. It asserts that any
contractible 2-complex 3-deforms to a point. Thus the above
theorem says that any contractible spine 3-deforms to a point. It is
not hard to see that a Heegaard diagram for $S^3$
corresponds to a contractible spine. Therefore the condition in
the above theorem is basically the same as in Theorem 3.2. It is
however well known (\cite{Wr}) that this geometric version of AC
conjecture is equivalent to what is actually referred to by many
as ``weak Andrews-Curtis conjecture" (see \cite{Mm}, for example).
It requires an additional operation, namely the addition (and deletion) of a new
generator and a new relator that is equal to the new generator. Thus Theorem 3.2 is stronger
than the theorem above. And our proof is quite different from the simple-homotopy theoretic
arguments sketched by Rolfsen.

Since the Poincare conjecture is true, as is now widely
believed after Perelman's work, we have

{\bf Corollary 3.3.} {\em If $(R, V)$ is a Heegaard diagram for a
homotopy 3-sphere and $\{a_{1}, \cdots, a_{n}\}$ is a set of
elements of $\pi_{1}(V, p)$ represented by $R=\{r_{1}, \cdots,
r_{n}\}$. Then $\{a_{1}, \cdots, a_{n}\}$ can be transformed into
a set of free generators of the free group $\pi_{1}(V, p)$ by a
finite number of Andrews-Curtis transformations. }

It is perhaps also of interest to see the implication of Theorem 3.2 and
Corollary 3.3 in pure algebraic terms. In the rest of this note, we
give an algebraic equivalence of a Heegaard diagram as defined in this note,
and then describe Corollary 3.3 in pure algebraic terms.

Let $F$ be a free group of rank $n$. A set $A=\{a_{1}, \cdots,
a_{k}\}$ of $k$ elements of $F$ is said to be conjugate to another
set $B=\{b_{1}, \cdots, b_{k}\}$ of $k$ elements of $F$ if for
each $i=1, \cdots, k$, $a_{i}$ is conjugate to $b_{i}$. A set
$\{c_{1}, \cdots, c_{n}\}$ of $n$ elements of $F$ is said to be
{\bf complete} or {\bf a complete set} if it is conjugate to a set
$\{a_{1}, \cdots, a_{n}\}$ for which there is a set
$\{b_{1},\cdots, b_{n}\}$ of elements of $F$ such that $F=(a_{1},
b_{1}, \cdots, a_{n}, b_{n})$ and $\prod_{i=1}^{n}[a_{i},
b_{i}]=1$.

We will now see that a complete set of $F$ is just the algebraic
equivalent of a complete system on the boundary of a handlebody of
genus $n$, that is,  a Heegaard diagram. Jaco \cite{Ja} proved that
every so-called splitting homomorphism is equivalent to a splitting
homomorphism induced by a Heegaard splitting of a closed 3-manifold.
This is used to prove
the following

{\bf Lemma 3.4. } {\em Let $F$ be a free group of rank $n$. A set
$\{c_{1}, \cdots, c_{n}\}$ of $n$ elements of $F$ is a complete
set if and only if there is an isomorphism $\alpha $ from $F$ to
$\pi_{1}(V, p)$ where $V$ is a handlebody and $p$ is a point on
$\partial V$ such that $\{\alpha(c_{1}), \cdots, \alpha(c_{n})\}$
is represented by a complete system on $\partial V$. }

{\bf Proof.} The ``if" part is straight forward. Without loss of
generality, assume that $F=\pi_{1}(V,p)$ for some handlebody $V$
of genus $n$ with a base point $p$ on $\partial V$ and $\{c_{1},
\cdots, c_{n}\}$ is represented by a complete system $R=\{r_{1},
\cdots, r_{n}\}$ on $\partial V$. Choose a complete system
$\tilde{R}=\{\tilde{r}_{1}, \cdots, \tilde{r}_{n}\}$ on $\partial
V$ such that $r_{i}\cap \tilde{r}_{j}$ is exactly one point if
$i=j$ and empty if $i\neq j$. Then $\pi_{1}(\partial V, p)=(A_{1},
B_{1}, \cdots, A_{n}, B_{n}\mid \prod_{i}[A_{i}, B_{i}])$ where
$A_{1}, \cdots, A_{n}$ are elements of $\pi_{1}(\partial V, p)$
determined respectively by $r_{1}, \cdots, r_{n}$ with some
connecting curves from $p$ and $B_{1}, \cdots, B_{n}$ by
$\tilde{r}_{1}, \cdots, \tilde{r}_{n}$. Let
$a_{i}=i_{V\ast}(A_{i})$ and $b_{i}=i_{V\ast}(B_{i})$, $i=1,
\cdots, n$, where $i_{V\ast}$ is the map induced by the inclusion
$i_{V}:\partial V\rightarrow V$. Then $F=\pi_{1}(V, p)=(a_{1},
b_{1},\cdots, a_{n}, b_{n})$ and $\prod_{i}[a_{i}, b_{i}]=1$.
Clearly $\{c_{1}, \cdots, c_{n}\}$ is conjugate to $\{a_{1},
\cdots, a_{n}\}$. This shows the ``if" part.

Now we show the ``only if" part. By the assumption, there are sets
$\{a_{1}, \cdots, a_{n}\}$ and $\{b_{1}, \cdots, b_{n}\}$ of
elements of $F$ such that $F=(a_{1}, b_{1}, \cdots, a_{n},
b_{n})$, $\prod_{i}[a_{i}, b_{i}]=1$ and $\{c_{1}, \cdots,
c_{n}\}$ is conjugate to $\{a_{1}, \cdots, a_{n}\}$. It suffices
to show that there is an isomorphism $\alpha $ from $F$ to
$\pi_{1}(V, p)$ where $V$ is a handlebody and $p$ is a point on
$\partial V$ such that $\{\alpha(a_{1}), \cdots, \alpha(a_{n})\}$
is represented by a complete system on $\partial V$.

Let $(S, p)$ be a closed surface of genus $n$ with base point $p$,
and $\{r_{1}, \cdots, r_{n}\}$ and $\{\tilde{r}_{1}, \cdots,
\tilde{r}_{n}\}$ be two complete systems on $S$ away from $p$ such
that $r_{i}\cap \tilde{r}_{j}$ is exactly one point if $i=j$ and
empty if $i\neq j$. Let $A_{1}, \cdots, A_{n}$ be elements of
$\pi_{1}(S, p)$ determined respectively by $r_{1}, \cdots, r_{n}$
with some connecting curves from $p$ and $B_{1}, \cdots, B_{n}$ be
elements determined by $\tilde{r}_{1}, \cdots, \tilde{r}_{n}$. We
may choose the connecting curves so that $\pi_{1}(S,
p)=(A_{1},B_{1} \cdots, A_{n}, B_{n}\mid \prod_{i}[A_{i},
B_{i}])$. Now let $G=(A_{1}, B_{1},\cdots, A_{n}, B_{n})$ be a
free group of rank $2n$ freely generated by the symbols $A_{1},
B_{1}, \cdots, A_{n}, B_{n}$ and $F_{1}=(y_{1}, \cdots, y_{n})$ be
a free group of rank $n$ freely generated by the symbols
$y_{1},\cdots, y_{n}$. Note that here we use each of the $A_{i}$'s
and $B_{i}$'s to denote two things: the element of $\pi_{1}(S,p)$
determined by one of the curves $r_{i}$ and $\tilde{r}_{i}$ with
some connecting curve, and a generating symbol for the free group
$G$. It should be clear what is meant from the context.

Let${\cal P}$ be the projection map from $G$ to $\pi_{1}(S, p)$,
${\cal P}(A_{i})=A_{i}$, ${\cal P}(B_{i})=B_{i}$, $i=1, \cdots,
n$. Define maps $\phi_{1}: G\rightarrow F_{1}$ and $\phi:
G\rightarrow F$ by $\phi_{1}(A_{i})=1,\phi_{1}(B_{i})=y_{i}$, and
$\phi(A_{i})=a_{i}, \phi(B_{i})=b_{i}$, $i=1, \cdots, n$. Both
$\phi_{1}$ and $\phi$ are surjective. Since $\phi_{1}(\prod_{i=1,
\cdots, n}[A_{i}, B_{i}])=\phi(\prod_{i=1, \cdots, n}[A_{i},
B_{i}])=1$, $\phi_{1}$ and $\phi$ both factor through the
fundamental group $\pi_{1}(S, p)$. So there are maps
$\psi_{1}:\pi_{1}(S, p)\rightarrow F_{1}$ and $\psi: \pi_{1}(S, p)
\rightarrow F$ such that $\phi_{1}=\psi_{1}\circ {\cal P}$ and $
\phi=\psi \circ {\cal P}$. Clearly $\psi_{1}$ and $\psi$ are also
surjective. Thus $\Psi=(\psi_{1}, \psi): \pi_{1}(S, p)\rightarrow
F_{1}\times F $ is a so-called splitting homomorphism (see
\cite{Ja}). By a theorem of Jaco \cite{Ja}, the map
$\Psi=(\psi_{1}, \psi): \pi_{1}(S, p)\rightarrow F_{1}\times F $
is equivalent to a splitting homomorphism $(i_{U\ast},
i_{V\ast}):\pi_{1}(T=U\cap V, q)\rightarrow \pi_{1}(U, q) \times
\pi_{1}(V, q) $ induced by a Heegaard splitting $(U, V)$ of some
closed 3-manifold $M$, where $q$ is some point on $T=U\cap V$.
This means that there are isomorphisms $\mu: \pi_{1}(S, p)
\rightarrow \pi_{1}(T, q)$, $\alpha_{1}: F_{1}\rightarrow
\pi_{1}(U, q)$ and $\alpha: F\rightarrow \pi_{1}(V, q)$ such that
$\alpha_{1}\circ \psi_{1}=i_{U\ast}\circ \mu $ and $\alpha \circ
\psi=i_{V\ast}\circ \mu$. By composing $\mu$ with an inner
automorphism of $\pi_{1}(T, q)$, and $\eta_{1}$ and $\eta$ with
the corresponding inner automorphisms of $\pi_{1}(U, q)$ and
$\pi_{1}(V, q)$ respectively if necessary, we can assume that
$\mu$ is induced by a homeomorphism $h: (S, p)\rightarrow (T, q)$.
Then $\alpha(a_{1}), \cdots, \alpha(a_{n})$ are elements of
$\pi_{1}(V, q)$ determined respectively by the simple closed
curves $h(r_{1}),\cdots, h(r_{n})$ on $T=\partial V$ with some
connecting curves. That is, $\alpha$ is an isomorphism from $F$ to
$\pi_{1}(V, p)$ such that $\{\alpha(a_{1}), \cdots,
\alpha(a_{n})\}$ is represented by the complete system
$\{h(r_{1}), \cdots, h(r_{n})\}$ on $\partial V$. This completes
the proof of the ``only if" part of Lemma 3.4 hence also that of
Lemma 3.4.

In practice, it seems not easy to determine whether or not a given
set of $n$ elements of $F$ is complete. We are most interested in
sets that generate $F$ normally. Following \cite{Bm}, we will call
a set $A=\{a_{1}, \cdots, a_{n}\}$ of $n$ elements of $F$ {\bf an
annihilating $n$-tuple} if $F=\langle a_{1}, \cdots,
a_{n}\rangle$.

{\bf Theorem 3.5.} {\em Let $F$ be a free group of rank $n$ and
$\{a_{1}, \cdots, a_{n}\}$ be an annihilating $n$-tuple for $F$.
If $\{a_{1}, \cdots, a_{n}\}$ is complete, then (assuming that the
Poincare conjecture is true) $\{a_{1}, \cdots, a_{n}\}$ can be
transformed into any free basis of $F$ by a finite number of
Andrews-Curtis transformations. }

{\bf Proof.} Suppose that $\{a_{1}, \cdots, a_{n}\}$ is complete.
By Lemma 3.4, there is an isomorphism $\alpha: F\mapsto \pi_{1}(V,
p)$ where $V$ is a handlebody of genus $n$ and $p$ is a point on
$\partial V$ such that $\{\alpha(a_{1}), \cdots, \alpha(a_{n})\}$
is represented by a complete system $R=\{r_{1}, \cdots, r_{n}\}$
on $\partial V$. Using $\alpha$ we identify $F$ with $\pi_{1}(V,
p)$ and $a_{1}, \cdots, a_{n}$ with elements of $\pi_{1}(V, p)$
represented respectively by $r_{1}, \cdots, r_{n}$.

Consider the Heegaard diagram $(R, V)$ where $R=\{r_{1}, \cdots,
r_{n}\}$. Since $\pi_{1}(V, p)=F=\langle a_{1}, \cdots,
a_{n}\rangle$, the manifold $M$ determined by $(R, V)$ is simply
connected. Theorem 3.5 then follows from Corollary 3.3.

\vspace{3cm}

Beijing, China

E-mail:  guoguangyuan@yahoo.com


\begin{thebibliography}{999}
\bibitem{DR} D. Rolfsen, The Poincare conjecture and its cousins, a survey note
         circulating on the internet.
\end{thebibliography}
\end{document}